\documentclass[11pt]{article}
\usepackage{amsmath,amsfonts,amssymb,amsthm,graphicx,verbatim,subfig,mathtools,wasysym}
\usepackage[all]{xy}
\usepackage{dsfont,enumitem}
\usepackage[a4paper,left=25mm]{geometry}
\usepackage[usenames]{color}
\usepackage{colortbl}

\ifx\pdftexversion\undefined

\usepackage[a4paper,colorlinks,linkcolor=black,filecolor=black,citecolor=black,urlcolor=black,pdfstartview=FitH]{hyperref}
\fi

\font\sixbb=msbm6
\font\eightbb=msbm8
\font\twelvebb=msbm10 scaled 1095
\newfam\bbfam
\textfont\bbfam=\twelvebb \scriptfont\bbfam=\eightbb
                           \scriptscriptfont\bbfam=\sixbb
\def\bb{\fam\bbfam\twelvebb}

\newcommand{\FF}{{\bb F}}

\newcommand{\bbb}{\mathbb{B}}
\newcommand{\bbg}{\mathbb{G}}

\newtheorem{theorem}{\bf Theorem}[section]
\newtheorem{claim}[theorem]{\bf Claim}

\newtheorem{proposition}[theorem]{\bf Proposition}
\newtheorem{corollary}[theorem]{\bf Corollary}
\newtheorem{remark}[theorem]{\bf Remark}
\newtheorem*{thw}{Theorem}
\newtheorem*{tha}{Theorem}
\newtheorem*{thb}{Theorem}
\newcommand{\enp}{\begin{flushright} $\Box$ \end{flushright}}
\newcommand{\beq}[0]{\begin{equation}}
\newcommand{\enq}[0]{\end{equation}}

\newcommand{\cb}{{\cal B}}

\newcommand{\rk}{{\rm rk}}

\newcommand{\cg}{\mathcal{G}}

\newcommand{\rrk}{\text{rk} \,}

\newcommand{\ccm}{\mathcal{M}}

\newcommand{\pf}{\text{Pf}}

\newcommand{\sgn}{\text{sgn}}
\newcommand{\prk}{{\rm prk}}
\newcommand{\per}{{\rm per}}
\newcommand{\spn}{{\rm span}}
\newcommand{\one}{\mathds{1}}
\newcommand{\tilq}{\tilde{q}}

\title{Maximal Generalized Rank in Graphical Matrix Spaces}
\author{Alexander Guterman \thanks{Department of Mathematics, Bar-Ilan University, Ramat-Gan 5290002, Israel. e-mail:  alexander.guter\-man@biu.ac.il}
\and Roy Meshulam \thanks{Department of Mathematics,
Technion, Haifa 32000, Israel. e-mail:
meshulam@technion.ac.il~. Supported by ISF grant 686/20.}
\and Igor Spiridonov  \thanks{Skolkovo Institute of Science and Technology, Skolkovo 121205, Russia; National Research University Higher School of Economics, Moscow 119048, Russia; Moscow Center for Fundamental and Applied Mathematics, Moscow 119991, Russia. e-mail: spiridonovia@yandex.ru}}

\begin{document}
\maketitle
\pagestyle{plain}

\begin{abstract}
Let $M_{n}(\FF)$ be the space of $n \times n$ matrices over a field $\FF$. For a subset $\cb \subset [n]^2$ let
$M_{\cb}(\FF)=\{A \in M_{n}(\FF):A(i,j)=0 {\rm ~for~} (i,j) \not\in \cb\}$. Let $\nu_b(\cb)$ denote the matching number of the $n$ by $n$ bipartite graph determined by $\cb$. For $S \subset M_n(\FF)$ let
$\rho(S)=\max\{\rk(A):A \in S\}$. Li, Qiao, Wigderson,
Wigderson and Zhang (arXiv:2206.04815, 2022) have recently proved the following characterization of the maximal dimension of bounded rank subspaces of $M_{\cb}(\FF)$.
\begin{thw}[Li, Qiao, Wigderson,
Wigderson, Zhang]
For any $\cb \subset [n]^2$
\begin{equation}
\label{e:abst1}
\max \left\{\dim W: ~W \leq M_{\cb}(\FF)~,~ \rho(W) \leq k\right\}=
\max \left\{|\cb'|:~ \cb' \subset \cb~,~ \nu_b(\cb') \leq k\right\}.
\end{equation}
\end{thw}
\noindent
The main results of this note are two extensions of (\ref{e:abst1}). Let $\mathbb{S}_n$ denote the symmetric group on $[n]$. For $\omega:\coprod_{n=1}^{\infty} \mathbb{S}_n \rightarrow \FF^*=\FF\setminus \{0\}$
define a function $D_{\omega}$ on each $M_n(\FF)$ by
$D_{\omega}(A)=\sum_{\sigma \in \mathbb{S}_n} \omega(\sigma) \prod_{i=1}^n A(i,\sigma(i))$.
Let $\rk_{\omega}(A)$ be the maximal $k$ such that there exists a $k \times k$ submatrix $B$ of $A$ with $D_{\omega}(B) \neq 0$. For $S \subset M_n(\FF)$
let $\rho_{\omega}(S)=\max \{\rk_{\omega}(A):A \in S\}$.
The first extension of (\ref{e:abst1}) concerns general weight functions.
\begin{tha}
For any $\omega$ as above and $\cb \subset [n]^2$
\begin{equation*}
\label{e:anyw2}
\max \left\{\dim W: ~W \leq M_{\cb}(\FF)~,~\rho_{\omega}(W) \leq k\right\}=
\max \left\{|\cb'|: \cb' \subset \cb~,~ \nu_b(\cb') \leq k\right\}.
\end{equation*}
\end{tha}
\noindent
Let $A_n(\FF)$ denote the space of alternating matrices in $M_n(\FF)$. For a graph $\cg \subset \binom{[n]}{2}$
let $A_{\cg}(\FF)=\{A \in A_{n}(\FF):A(i,j)=0 {\rm ~if~} \{i,j\} \not\in \cg\}$.
Let $\nu(\cg)$ denote the matching number of $\cg$.
The second extension of (\ref{e:abst1}) concerns general graphs.
\begin{thb}
For any $\cg \subset \binom{[n]}{2}$
\begin{equation*}
\label{e:abst3}
\max \left\{\dim U: ~U \leq A_{\cg}(\FF)~,~ \rho(U) \leq 2k\right\}=
\max \left\{|\cg'|:~ \cg' \subset \cg~,~ \nu(\cg') \leq k\right\}.
\end{equation*}
\end{thb}

\end{abstract}

\section{Introduction}
\label{s:intro}

Let $M_{n}(\FF)$ denote the space of $n \times n$ over a field $\FF$. For $A \in M_n(\FF)$ and
subsets $\emptyset \neq I=\{i_1,\ldots,i_k\}, J=\{j_1,\ldots,j_k\} \subset [n]:=\{1,\ldots,n\}$ 
such that $i_1<\cdots < i_k$
and $j_1 < \cdots < j_k$ let $B=A[I|J] \in M_k(\FF)$ be given by
$B(\alpha,\beta)=A(i_{\alpha},j_{\beta})$ for
$1 \leq \alpha,\beta \leq k$. For a vector space $W$ we write $V \leq W$ if $V$ is a linear subspace of $W$. 
Let $\binom{[n]}{k}$ denote the family of all $k$-element subsets of $[n]$.
For column vectors $u_1,u_2 \in \FF^n$ let $u_1 \otimes u_2=u_1\cdot u_2^t \in M_n(\FF)$.
The tensor product of two linear subspaces $U_1,U_2 \leq \FF^n$ is given by
$
U_1\otimes U_2=
\spn\left\{u_1 \otimes u_2: u_1 \in U_1, u_2 \in U_2\right\}.
$
For a subset $S \subset M_n(\FF)$ let $\rho(S)=\max\{\rrk(A): A \in S\}$ denote the maximal rank of a matrix in $S$. The following result was proved by Flanders \cite{Flanders62} under the assumption $|\FF| \geq k+1$, and in \cite{M85} for all fields.
\begin{theorem}[\cite{Flanders62,M85}]
\label{t:fm}
Let $W \leq M_n(\FF)$ be a linear subspace such that $\rho(W) \leq k$. Then:
(i) $\dim W \leq kn$. (ii) $\dim W=kn$ iff $W=U \otimes \FF^n$ or $W=\FF^n \otimes U$ for some $k$-dimensional linear subspace
$U \leq \FF^n$.
\end{theorem}
\noindent
For $i \in [n]$ let $e_i$ denote the $i$-th unit vector in $\FF^n$.
For a subset $\cb \subset [n]^2$ let
\begin{equation*}
\label{e:mcb}
M_{\cb}(\FF)=\spn\left\{e_i \otimes e_j: (i,j) \in \cb\right\}.
\end{equation*}
\noindent
A \emph{bipartite matching} in $\cb$ is a subset $\cb_0 \subset \cb$ such that if
$(i,j) \neq (i',j') \in \cb_0$ then $i\neq i'$ and $j \neq j'$.
The \emph{bipartite matching number} of $\cb$ is
\[
\nu_b(\cb):=\max \left\{|\cb_0|: \cb_0 {\rm ~is~a~bipartite~matching~in~} \cb\right\}.
\]
\noindent
Li, Qiao, Wigderson, Wigderson and Zhang \cite{LQWWZ22} have recently established the following
\begin{theorem}[\cite{LQWWZ22}]
\label{t:lqwwz}
For any $\cb \subset [n]^2$
\begin{equation}
\label{e:lqwwz}
\max \left\{\dim W: ~W \leq M_{\cb}(\FF)~,~ \rho(W) \leq k\right\}=
\max \left\{|\cb'|:~ \cb' \subset \cb~,~ \nu_b(\cb') \leq k\right\}.
\end{equation}
\end{theorem}
\begin{remark}
\label{r:1}
When $\cb=[n]^2$ Theorem \ref{t:lqwwz} specializes to Theorem \ref{t:fm}(i). Indeed, K\H{o}nig's Theorem
(see e.g. Theorem 3.1.11 in \cite{West96})
implies that $\max \left\{|\cb'|: \cb' \subset [n]^2, \nu_b(\cb') \leq k\right\}=kn$.
\end{remark}
\noindent
In this paper we give some extensions of Theorems \ref{t:fm} and \ref{t:lqwwz}.
Let $\mathbb{S}_n$ denote the symmetric group on $[n]$. For a weight
function $\omega:\coprod_{n=1}^{\infty} \mathbb{S}_n \rightarrow \FF^*=\FF\setminus \{0\}$
let $D_{\omega}:\coprod_{n=1}^{\infty} M_n(\FF) \rightarrow \FF$ be defined on $A=(A(i,j))_{i,j=1}^n \in M_n(\FF)$ by
\begin{equation*}
\label{e:dalpha}
D_{\omega}(A)=\sum_{\sigma \in \mathbb{S}_n} \omega(\sigma) \prod_{i=1}^n A(i,\sigma(i)).
\end{equation*}
\begin{remark}
The functions $D_{\omega}$ were considered by de Seguins Pazzis under the name of
\emph{Schur matrix functionals}. See his paper \cite{DSP19} for an in-depth study of the linear preservers of
various $D_{\omega}$'s.
\end{remark}
\noindent
The {\it $\omega$-rank} $\rk_{\omega}(A)$ of a matrix $A \in M_n(\FF)$ is the maximal $k$ such that there exist
$I,J \in \binom{[n]}{k}$ such that
$D_{\omega}\left(A[I|J]\right) \neq 0$. For $S \subset M_n(\FF)$
let $\rho_{\omega}(S)=\max \{\rk_{\omega}(A):A \in S\}$.
Note that for the sign function
$D_{\sgn}(A)=\det A$ and $\rk_{\sgn}(A)=\rk(A)$. Let $\one$ be the constant function $\one(\sigma)\equiv 1$.
Then $D_{\one}(A)= \per A$ is the permanent of $A$, and $\prk(A):=\rk_{\one}(A)$ is the {\it permanental rank} of $A$.  See \cite{GutSp} for applications of permanental rank to linear preserver problems. Our first result is an extension of Theorem \ref{t:lqwwz} to general weight functions.
\begin{theorem}
\label{t:anyw}
For any $\omega:\coprod_{n=1}^{\infty} \mathbb{S}_n \rightarrow \FF^*$ and any $\cb \subset [n]^2$
\begin{equation}
\label{e:anyw}
\max \left\{\dim W: ~W \leq M_{\cb}(\FF)~,~\rho_{\omega}(W) \leq k\right\}=
\max \left\{|\cb'|: \cb' \subset \cb~,~ \nu_b(\cb') \leq k\right\}.
\end{equation}
\end{theorem}
\noindent
Specializing to $\cb=[n]^2$ and $\omega=\one$ we obtain the following permanental counterpart of Theorem
\ref{t:fm}(i).
\begin{corollary}
\label{c:perm}
Let $W \leq M_n(\FF)$ be a linear subspace such that $\prk(A) \leq k$ for all $A \in W$.
Then $\dim W \leq kn$.
\end{corollary}
\noindent
For the fields of characteristic $2$, we have $\det=\per$ and thus the equality cases in Corollary \ref{c:perm} are those given in Theorem \ref{t:fm}(ii).
For $(k,n)=(1,2)$ it can be checked that the only $2$-dimensional subspaces $W \leq M_2(\FF)$ such that
$\dim W=2$ and $\rho_{\one}(W)=1$ are $W_u$ and its transpose, where for $0 \neq u=(a,b) \in \FF^2$
\[
W_u=\left\{\left(
\begin{array}{rr}
ax & ay \\
-bx  & by
\end{array}
\right) : x,y \in \FF\right\}.
\]
\noindent
In general we have the following
\begin{theorem}
\label{t:eqcase}
Suppose that ${\rm char}\, \FF \neq 2$, $k \leq n$ and $n \geq 3$. Then $W \leq M_n(\FF)$ satisfies $\rho_{\one}(W)=k$ and $\dim W=kn$ iff $W=\spn\{e_i\}_{i \in I} \otimes \FF^n$ or
$W=\FF^n \otimes \spn\{e_i\}_{i \in I}$ for some $I \in \binom{[n]}{k}$.
\end{theorem}
A matrix $A=\big(A(i,j)\big)_{i,j=1}^n \in M_n(\FF)$
is \emph{alternating} if $A=-A^t$ and $A(i,i)=0$ for $1 \leq i \leq n$. Let $A_n(\FF)$ denote the space of alternating matrices in $M_n(\FF)$. Recall that $\rk(A)$ is even for all $A \in A_n(\FF)$, and that
$\rk(A)=2k$ iff there exists a subset $I \in \binom{[n]}{2k}$ such that the principal submatrix
$A[I|I]$  is nonsingular.
For $u,v \in \FF^n$ let $u \wedge v= u \otimes v -v \otimes u \in A_n(\FF)$.
For a subset $\cg \subset \binom{[n]}{2}$ let
\[
A_{\cg}(\FF)=\spn\left\{e_i \wedge e_j: \{i,j\} \in \cg\right\}.
\]
\noindent
A \emph{matching} in $\cg$ is a subset $\cg_0 \subset \cg$ such that $f \cap f' =\emptyset$ for all
$f \neq f'\in \cg_0$.
The \emph{matching number} of $\cg$ is
\[
\nu(\cg):=\max \left\{|\cg_0|: \cg_0 {\rm ~is~a~matching~in~} \cg \right\}.
\]
Our final result is an extension of Theorem \ref{t:lqwwz} to spaces of alternating
matrices supported on general graphs.
\begin{theorem}
\label{t:anygr}
For any $\cg \subset K_n:=\binom{[n]}{2}$
\begin{equation}
\label{e:anygr}
\max \left\{\dim U: ~U \leq A_{\cg}(\FF)~,~ \rho(U) \leq 2k\right\}=
\max \left\{|\cg'|:~ \cg' \subset \cg~,~ \nu(\cg') \leq k\right\}.
\end{equation}
\end{theorem}
\begin{remark}
\label{r:gimplyb}
Theorem \ref{t:anygr} implies Theorem \ref{t:lqwwz}. Indeed, given $\cb \subset [n]^2$ let
\[
\cg=\left\{ \{i,j+n\}: (i,j) \in \cb\right\}\subset \binom{[2n]}{2}.
\]
It is straightforward to check that
\[
\max \left\{\dim W: ~W \leq M_{\cb}(\FF)~,~ \rho(W) \leq k\right\}=
\max \left\{\dim U: ~U \leq A_{\cg}(\FF)~,~ \rho(U) \leq 2k\right\}
\]
and
\[
\max \left\{|\cb'|:~ \cb' \subset \cb~,~ \nu_b(\cb') \leq k\right\}=
\max \left\{|\cg'|:~ \cg' \subset \cg~,~ \nu(\cg') \leq k\right\}.
\]
\noindent
Thus (\ref{e:lqwwz}) follows from (\ref{e:anygr}).
\end{remark}

The paper is organized as follows. In Section \ref{s:bndw} we prove Theorem \ref{t:anyw}. Our main tool is a combinatorial lower bound on the maximal $\omega$-rank in a subspace of matrices given in Proposition \ref{p:mainp}.
In Section \ref{s:maxprk} we prove Theorem \ref{t:eqcase}. In Section \ref{s:mrgsa} we use a result from \cite{M17} to establish Theorem \ref{t:anygr}. We conclude in Section \ref{s:conc} with
some remarks and open problems.

\section{Maximal $\omega$-Rank in Subspaces of $M_{\cb}(\FF)$}
\label{s:bndw}

In this section we prove Theorem \ref{t:anyw}. We will need the following facts.
\begin{claim}
\label{c:rwmb}
For any $\cb \subset [n]^2$
\begin{equation}
\label{e:rwmb}
\rho_{\omega}(M_\cb(\FF))=\nu_b(\cb).
\end{equation}
\end{claim}
\noindent
{\bf Proof.} Let $\rho_{\omega}(M_\cb(\FF))=k$ and $\nu_b(\cb)=\ell$.
Then there exists a matrix
$A \in M_\cb(\FF)$ and
a $k \times k$ submatrix $A'=A[I|J]$ such that $D_{\omega}(A') \neq 0$. Let $I=\{i_1< \cdots <i_k\}$
and $J=\{j_1< \cdots <j_k\}$. Then
\begin{equation*}
0 \neq D_{\omega}(A') =\sum_{\pi \in \mathbb{S}_k} \omega(\pi) \prod_{t=1}^k A(i_t, j_{\pi(t)}).
\end{equation*}
It follows that there exists $\pi \in \mathbb{S}_k$ such that $A(i_t, j_{\pi(t)}) \neq 0$ for
all $1 \leq t \leq k$.
Thus $\cb_0=\left\{ (i_t,j_{\pi(t)})\right\}_{t=1}^k$
is a bipartite matching of size $k$ in $\cb$ and therefore $\ell \geq k$. For the other direction, let
$\cb_0=\left\{(i_t,j_t)\right\}_{t=1}^{\ell}$ be a bipartite matching of size $\ell$ in $\cb$.
By reordering we may assume that $i_1< \cdots <i_{\ell}$. Let $\pi \in S_{\ell}$ such that
$j_{\pi(1)} < \cdots < j_{\pi(\ell)}$. Let $I=\{i_1,\ldots,i_{\ell}\}$,
$J=\{j_1,\ldots,j_{\ell}\}$, and let $A=\sum_{t=1}^{\ell} e_{i_t} \otimes e_{j_t} \in M_\cb(\FF)$.
Then $A'=A[I|J]$ satisfies $D_{\omega}(A')=\omega(\pi^{-1}) \neq 0$ and therefore
$k \geq \rk_{\omega}(A)=\ell$.
{\enp}
\noindent
Let $\prec$ be the lexicographic order on $[n]^2$, i.e. $(i,j)\prec (i',j')$ if
either $i<i'$ or $i=i'$ and $j<j'$. For $0 \neq A \in M_n(\FF)$ let $q(A)=\min_{\prec} \{(i,j): A(i,j) \neq 0\}$. For $S \subset M_n(\FF)$ let $\bbb(S)=\left\{q(A): A \in S\right\}$. Note that if $W \leq M_n(\FF)$ is a linear subspace then $|\bbb(W)|=\dim W$.
The main tool in the proof of Theorem \ref{t:anyw} is the following
\begin{proposition}
\label{p:mainp}
Let $W \leq M_n(\FF)$ be a linear subspace. Then
\begin{equation}
\label{e:mainp}
\rho_{\omega}(W) \geq \nu_b(\bbb(W)).
\end{equation}
\end{proposition}
\begin{remark}
\label{r:2}
The case $\omega=\sgn$ of Proposition \ref{p:mainp} is equivalent to Theorem 1 in \cite{M85}. The proof for general weight functions given below requires an additional idea.
\end{remark}
\noindent
{\bf Proof of Proposition \ref{p:mainp}.} Let $k=\nu_b(\bbb(W))$.
Then there exist $A_1,\ldots,A_k \in W$ with $q(A_t)=(i_t,j_t)$ such that
$\{(i_t,j_t)\}_{t=1}^k$ is a bipartite matching of size $k$.
Let $I=\{i_1,\ldots,i_k\}, J=\{j_1,\ldots,j_k\}$.
By reordering and rescaling the matrices $A_t$'s, we may assume that $i_1< \cdots < i_{k}$ and that
 $A_t(i_t,j_t)=1$. Let $\pi \in \mathbb{S}_{k}$
such that $j_{\pi(1)}<\cdots  <j_{\pi(k)}$. For $1 \leq t \leq k$ let
$C_t=A_{t}[I|J] \in M_{k}(\FF)$. Note that
\begin{equation}
\label{e:cab}
C_t(\alpha,\beta)=\left\{
\begin{array}{ll}
0 & \alpha<t ~{\rm or}~ (\alpha=t ~\&~ \beta<\pi^{-1}(t)), \\
1& (\alpha,\beta)=(t,\pi^{-1}(t)).
\end{array}
\right.
\end{equation}
Let $x=(x_1,\ldots,x_{k})$.  Let
 $G(x)=\sum_{t=1}^{k}x_t C_t$ and consider the polynomial
\begin{equation*}
g(x)=D_{\omega} \left(G(x)\right)=\sum_{\sigma \in \mathbb{S}_{k}} \omega(\sigma) \prod_{\ell=1}^{k}G(x)(\ell,\sigma(\ell))
 \in \FF[x_1,\ldots,x_k].
\end{equation*}
\begin{claim}
\label{c:prec}
There exists a $\lambda \in \FF^{k}$
such that $g(\lambda) \neq 0$.
\end{claim}
\noindent
We will use Alon's Combinatorial Nullstellensatz (Theorem 1.2 in \cite{Alon99}).
\begin{theorem}[Alon \cite{Alon99}]
\label{t:null}
Let $\FF$ be an be an arbitrary field and let $g=g(x_1,\ldots,x_k) \in \FF[x_1\ldots,x_k]$.
Suppose the total degree $\deg(g)$ of $g$ is $\sum_{t=1}^k d_t$ where each $d_t$ is a nonnegative
integer, and suppose the coefficient of $\prod_{t=1}^k x_t^{d_t}$ in $g$ is nonzero.
Then, if $\Lambda_1,\ldots,\Lambda_k$ are subsets of $\FF$ with $|\Lambda_t|>d_t$, there exist $\lambda_1 \in \Lambda_1,\ldots, \lambda_k \in \Lambda_t$ such that
$g(\lambda_1,\ldots,\lambda_k) \neq 0$.
\end{theorem}
\noindent
{\bf Proof of Claim \ref{c:prec}:} We will show that
the monomial $x_1 \cdots x_{k}$ appears with a nonzero coefficient in $g(x)$.
Indeed, let $\sigma \in \mathbb{S}_{k}$. By (\ref{e:cab}), for any $1 \leq t \leq k$ the variable
$x_t$ does not appear in $\prod_{\ell <t} G(x)\left(\ell,\sigma(\ell)\right)$.
It follows that the coefficient of
 $x_1 \cdots x_{k}$ in $\prod_{\ell=1}^{k} G(x)\left(\ell,\sigma(\ell)\right)$ is
 \begin{equation*}
\gamma(\sigma):= \prod_{\ell=1}^{k} C_{\ell}(\ell,\sigma(\ell)).
 \end{equation*}
 \noindent
Let $1 \leq \ell \leq k$.
If $\sigma(\ell)<\pi^{-1}(\ell)$ then $C_{\ell}(\ell,\sigma(\ell))=0$ by (\ref{e:cab}). It follows that
if $\gamma(\sigma) \neq 0$ then $\sigma(\ell) \geq \pi^{-1}(\ell)$ for all $1 \leq \ell \leq k$, i.e.
$\sigma=\pi^{-1}$. Therefore the coefficient of $x_1 \cdots x_{k}$ in $g(x)$ is
\begin{equation*}
\begin{split}
\omega(\pi^{-1})\gamma(\pi^{-1})= \omega(\pi^{-1})\prod_{\ell=1}^{k} C_{\ell}(\ell,\pi^{-1}(\ell))=\omega(\pi^{-1}) \neq 0.
\end{split}
 \end{equation*}
Applying Theorem \ref{t:null} for the polynomial $g$, with $d_1=\cdots=d_{k}=1$ and
$\Lambda_1=\cdots=\Lambda_{k}=\{0,1\}$, it follows that there exists $\lambda \in \{0,1\}^{k}$
such that $D_{\omega}\left(G(\lambda)\right)=g(\lambda) \neq 0$. As $G(\lambda)$ is a $k \times k$
submatrix of $\sum_{t=1}^{k} \lambda_t A_{t} \in W$  it follows that
$\rho_{\omega}(W)\geq k$.
{\enp}
\noindent
{\bf Proof of Theorem \ref{t:anyw}.} The $\geq$ direction of (\ref{e:anyw}) follows
from Claim \ref{c:rwmb}. Indeed,
if $\cb' \subset \cb$ satisfies $\nu_b(\cb') \leq k$, then $W=M_{\cb'}(\FF) \leq M_{\cb}(\FF)$ satisfies $\rho_{\omega}(W)=\nu_b(\cb') \leq k$ and $\dim W=|\cb'|$. For the $\leq$ direction, suppose
$W \leq M_{\cb}(\FF)$ satisfies $\rho_{\omega}(W) \leq k$. Let $\cb'=\bbb(W) \subset \cb$. Then
$|\cb'|=\dim W$ and Proposition \ref{p:mainp} implies that $\nu_b(\cb') \leq \rho_{\omega}(W) \leq k$.
{\enp}

\section{Maximal Dimensional Spaces of Bounded Permanental Rank}
\label{s:maxprk}
In this section we prove Theorem \ref{t:eqcase}. The main ingredient of the argument is the following
\begin{proposition}
\label{p:kkpo}
Let $k \geq 2$ and suppose that ${\rm char}\, \FF \neq 2$
and $U \leq M_{k+1}(\FF)$ satisfies $\bbb(U)=[k] \times [k+1]$ and
$\rho_{\one}(U)=k$. Then $U=\spn\{e_i\}_{i \in [k]} \otimes \FF^{k+1}$.
\end{proposition}
\noindent
{\bf Proof.} The assumption $\bbb(U)=[k] \times [k+1]$ implies that $U$ has a basis
\[
\left\{A_{ij}: (i,j) \in [k] \times [k+1]\right\}
\]
such that
\[
A_{ij}=e_i \otimes e_j
+\sum_{\ell=1}^{k+1} \lambda_{ij\ell} e_{k+1} \otimes e_{\ell}
\]
for some $\lambda_{ij\ell} \in \FF$.
\begin{claim}
\label{c:mzr}
$\lambda_{ij\ell}=0$ for any $1 \leq i \leq k$ and $1 \leq j \neq \ell \leq k+1$.
\end{claim}
\noindent
{\bf Proof.} Suppose for contradiction that
$\lambda_{ij\ell} \neq 0$ for some $1 \leq i \leq k$ and $1 \leq j \neq \ell \leq k+1$.
Let
\[
[k] \setminus \{i\}=\{i_1,\ldots,i_{k-1}\}~~,~~
[k+1] \setminus \{\ell,j\}=\{j_1,\ldots,j_{k-1}\},
\]
\noindent
and for $\theta \in \FF$ let
\[
C_{\theta}=\theta A_{ij}+\sum_{t=1}^{k-1}A_{i_t j_t} \in U.
\]
\noindent
Clearly, the only permutation $\pi \in S_{k+1}$ that corresponds to a nonzero term in the expansion of
$\per\,C_{\theta}$ is that given by $\pi(i)=j,\pi(k+1)=\ell$ and $\pi(i_t)=j_t$ for $1 \leq t \leq k-1$.
Therefore
\begin{equation}
\label{e:mzr}
\begin{split}
\per\, C_{\theta}&=
C_{\theta}(i,j)\cdot C_{\theta}(k+1,\ell) \cdot \prod_{t=1}^{k-1} C_{\theta}(i_t,j_t) \\
&=C_{\theta}(i,j)\cdot C_{\theta}(k+1,\ell)  \\
&=\theta \left(\theta A_{ij}(k+1,\ell)+\sum_{t=1}^{k-1}A_{i_t j_t}(k+1,\ell)\right) \\
&=\theta \left(\theta \lambda_{ij\ell}+\sum_{t=1}^{k-1} \lambda_{i_tj_t \ell}\right).
\end{split}
\end{equation}
Eq. (\ref{e:mzr}) and the assumption $\lambda_{ij\ell} \neq 0$ imply that $\per \, C_{\theta}$ is a nonzero polynomial of degree $2$ in $\theta$. As $|\FF| \geq 3$ there exists $\theta \in \FF$ such that
$\per \, C_{\theta} \neq 0$, contradicting $\rho_{\one}(U)=k$.
{\enp}
\noindent
For $(i,j) \in [k] \times [k+1]$ let $\mu_{ij}=\lambda_{ijj}$. Claim \ref{c:mzr}
implies that $A_{ij}=e_i \otimes e_j + \mu_{ij} e_{k+1} \otimes e_j$.
\begin{claim}
\label{c:ij}
$\mu_{ij}=0$ for all $(i,j) \in [k] \times [k+1]$.
\end{claim}
\noindent
{\bf Proof.} Fix $i \in [k]$ and $j' \neq j''  \in [k+1]$.
Let
\[
[k]\setminus \{i\}=\{i_1,\ldots,i_{k-1}\}~~,~~
[k+1] \setminus \{j',j''\}=\{j_1,\ldots,j_{k-1}\},
\]
and let
\[
C=A_{ij'}+A_{ij''}+ \sum_{t=1}^{k-1} A_{i_tj_t}.
\]
\noindent
Clearly, if $\pi \in S_{k+1}$ corresponds to a nonzero term in the expansion of
$\per\,C$ then $\pi(i_t)=j_t$ for $1 \leq t \leq k-1$ and either
$(\pi(i),\pi(k+1))=(j',j'')$ or $(\pi(i),\pi(k+1))=(j'',j')$.
It follows that
\begin{equation*}
\label{e:ijel}
\begin{split}
\per\, C&=
\left(C(i,j')\cdot C(k+1,j'')+C(i,j'')\cdot C(k+1,j')\right)
 \cdot \prod_{t=1}^{k-1} C(i_t,j_t) \\
&=C(k+1,j'')+C(k+1,j')=\mu_{ij''}+\mu_{ij'}.
\end{split}
\end{equation*}
\noindent
Together with the assumption that $\rho_{\one}(U)=k$ this implies that
\begin{equation}
\label{e:muij}
\mu_{ij''}+\mu_{ij'}=\per\, C=0.
\end{equation}
\noindent
As (\ref{e:muij}) holds for all $1 \leq j' \neq j'' \leq k+1$ and $k \geq 2$, the assumption
${\rm char}(\FF) \neq 2$ implies that
$\mu_{ij}=0$ for all $(i,j) \in [k] \times [k+1]$, thereby completing the proof of Claim \ref{c:mzr}
and of Proposition \ref{p:kkpo}.
{\enp}
\noindent
{\bf Proof of Theorem \ref{t:eqcase}.}
We may assume that $k<n$. Let $W \leq M_n(\FF)$ satisfy
$\rho_{\one}(W)=k$ and $\dim W=kn$. Proposition \ref{p:mainp} implies that
\[
\nu(\bbb(W)) \leq \rho_{\one}(W)=k.
\]
As $|\bbb(W)|=\dim W=kn$, it follows by K\H{o}nig's theorem that $\bbb(W)$ is either
$I \times [n]$ or $[n] \times I$ for some $I \in \binom{[n]}{k}$.
Consider the first case $\bbb(W)=I \times [n]$ and let $I=\{i_1,\ldots,i_k\}$ where
$1 \leq i_1< \cdots < i_k \leq n$. We have to show that $W=\spn\{e_i: i \in I\} \otimes \FF^n$.
Fix a pair $(i',j') \in ([n] \setminus I) \times [n]$. Define $i_{k+1}=i'$ and choose $J=\{j_1,\ldots,j_{k+1}\} \subset [n]$ such that $j_1<\cdots <j_{k+1}$ and $j'=j_{\beta'}$ for some $1 \leq \beta' \leq k+1$. For $A \in W$ let $\tilde{A} \in M_{k+1}(\FF)$ be given by $\tilde{A}(\alpha,\beta)=A(i_{\alpha},j_{\beta})$ for $(\alpha,\beta) \in [k+1]^2$.
Let $U=\{\tilde{A}:A \in W\}$. Clearly
$\bbb(U)=[k] \times [k+1]$ and
$\rho_{\one}(U)=k$. Proposition \ref{p:kkpo} then implies that $U=\spn\{e_{\alpha}\}_{\alpha \in [k]} \otimes \FF^{k+1}$. In particular $A(i',j')=\tilde{A}(k+1,\beta')=0$ for all $A \in W$. Therefore
\begin{equation}
\label{e:wspn}
W \leq \spn\{e_i\}_{i \in I} \otimes \FF^n.
\end{equation}
\noindent
As $\dim W=kn$ it follows that there is equality in (\ref{e:wspn}).
The case $\bbb(W)=[n] \times I$ for some $|I|=k$ is handled similarly.
{\enp}

\section{Maximal Rank in Subspaces of $A_{\cg}(\FF)$}
\label{s:mrgsa}
In this section we prove Theorem \ref{t:anygr}.
We first recall the definition of the {\it Pfaffian} of an alternating matrix of even order
$A \in A_{2k}(\FF)$.
A \emph{perfect matching} in $K_{2k}$ is a matching of size $k$.
Let $\ccm_{2k}$ denote the set of all perfect matchings in $K_{2k}$. For $M=\{f_1,\ldots,f_{k}\} \in \ccm_{2k}$
such that $f_t=\{i_t<j_t\}$ for $1 \leq t \leq k$ and $i_1<\cdots<i_k$, let
$$\theta(M)=\text{sgn}  \left(
\begin{array}{ccccc}
1 & 2 & \cdots & 2k-1 & 2k  \\
i_1 & j_1 & \cdots & i_k & j_k
\end{array}
\right)
$$
and let
$$\mu(A,M)= \prod_{t=1}^{k} A(i_t,j_t).$$
The Pfaffian of $A$ is defined by
\[
\pf(A)=\sum_{M \in \ccm_{2k}} \theta(M) \mu(A,M).
\]
It is well known that $\det(A)=\pf(A)^2$ (see e.g. Exercise 4.24 in \cite{Lovasz}).
We will need the following facts.
\begin{claim}
\label{c:rmba}
For any $\cg \subset \binom{[n]}{2}$
\begin{equation}
\label{e:rmba}
\rho(A_\cg(\FF))=2\nu(\cg).
\end{equation}
\end{claim}
\noindent
{\bf Proof.} Let $\rho(A_{\cg}(\FF)=2k$ and $\nu(\cg)=\ell$.
Let $\left\{\{i_1,j_1\},\ldots,\{i_{\ell},j_{\ell}\} \right\}$ be a matching in $\cg$. Then
$A=\sum_{t=1}^{\ell} e_{i_t} \wedge e_{j_t} \in A_{\cg}(\FF))$ and
$\rho(A) =2 \ell$. Therefore $k \geq \ell$. For the other direction, let
$A \in A_{\cg}(\FF)$ such that $\rk(A)=2k$. Then there exists $I=\{\alpha_1< \cdots <\alpha_{2k}\}$
such that $B=A[I|I] \in A_{2k}(\FF)$ is nonsingular. Thus
\[
0 \neq \det(B)=\pf(B)^2=\sum_{M \in \ccm_{2k}} \theta(M) \mu(B,M).
\]
\noindent
Hence there exists a matching $M=\{f_1,\ldots,f_k\} \in \ccm_{2k}$ such that $\mu(B,M) \neq 0$.
Writing $f_t=\{i_t<j_t\}$ for $1 \leq t \leq k$ it follows that
\[
0 \neq \mu(B,M)=\prod_{t=1}^{k} B(i_t,j_t)=\prod_{t=1}^{k} A(\alpha_{i_t},\alpha_{j_t}).
\]
Thus $\left\{\{\alpha_{i_t},\alpha_{j_t}\}: 1 \leq t \leq k\right\}$ is a matching of size $k$ in $\cg$
and therefore $\ell=\nu(\cg) \geq k$.
{\enp}
\noindent
For $A \in A_n(\FF)$ with $q(A)=(i,j)$ define $\tilq(A)=\{i,j\}$.
For $S \subset A_n(\FF)$ let $\bbg(S)=\left\{\tilq(A): A \in S\right\}$. Note that if $U \leq A_n(\FF)$ is a linear subspace then $|\bbg(U)|=\dim U$.
The key ingredient in the proof of Theorem \ref{t:anygr} is the following
result (Theorem 1.2 in \cite{M17}).
\begin{theorem}[\cite{M17}]
\label{t:alter}
Let $U \leq A_n(\FF)$ be a linear subspace.
Then $\rho(U) \geq 2\nu(\bbg(U))$.
\end{theorem}
\noindent
{\bf Proof of Theorem \ref{t:anygr}.} The $\geq$ direction of (\ref{e:anygr}) follows
from Claim \ref{c:rmba}. Indeed,
if $\cg' \subset \cg$ satisfies $\nu(\cg') \leq k$, then $U=A_{\cg'}(\FF) \leq A_{\cg}(\FF)$ satisfies $\rho(U)=2\nu(\cg') \leq 2k$ and $\dim U=|\cg'|$. For the $\leq$ direction, suppose
$U \leq A_{\cg}(\FF)$ satisfies $\rho(U) \leq 2k$. Let $\cg'=\bbg(U) \subset \cg$. Then
$|\cg'|=\dim U$ and Theorem \ref{t:alter} implies that $\nu(\cg') \leq \frac{\rho(U)}{2} \leq k$.
{\enp}

\section{Concluding Remarks}
\label{s:conc}
In this note we proved two extensions of the combinatorial characterization due to Li, Qiao, Wigderson, Wigderson and Zhang \cite{LQWWZ22} of the maximal dimension of bounded rank subspaces of the graphical matrix space $M_{\cb}(\FF)$ associated with a bipartite graph $\cb$. Theorem \ref{t:anyw} shows that the above characterization remains valid for a wide class of generalized rank functions, including e.g. the permanental rank. In a different direction, Theorem \ref{t:anygr} extends the characterization to bounded rank subspaces of the graphical matrix space $A_{\cg}(\FF)$ associated with a general graph $\cg$.
We conclude with the following two remarks.
\begin{itemize}
\item
Theorem \ref{t:fm} provides a classification of the spaces $W \leq M_n(\FF)$ such that
$\rho(W)=k$ and $\dim W=kn$. The analogous (but different) classification of
spaces $W \leq M_n(\FF)$ such that $\rho_{\one}(W)=k$ and $\dim W=kn$ is given in Theorem \ref{t:eqcase}.
It would be interesting to extend these results to bounded rank subspaces of $M_{\cb}(\FF)$ and
of $A_{\cg}(\FF)$ for various $\cb \subset [n]^2$ and $\cg \subset \binom{[n]}{2}$.

\item
For an element $u$ in the $p$-th exterior power $\bigwedge^p \FF^n$,
 let $E(u)$ denote the minimal subspace $U \leq \FF^n$ such that
$u \in \bigwedge^p U$. The {\it rank} of $u$ is $\rk(u)=\dim E(u)$. For sufficiently large fields,
Theorem \ref{t:alter} is the case $p=2$ of Theorem 2.1 in \cite{GM01} that gives a lower bound on
$\rho(U)$ for $U \leq \bigwedge^p \FF^n$ in terms of the weak matching number of a certain $p$-uniform hypergraph associated to $U$. It seems likely that this lower bound may be useful in obtaining extensions
of Theorem \ref{t:anygr} to structured subspaces of $\bigwedge^p \FF^n$ for general $p$.

\end{itemize}


\begin{thebibliography}{99}

\bibitem{Alon99}
N. Alon,
Combinatorial Nullstellensatz,
{\it Combin. Probab. Comput.}, {\bf 8}(1999) 7--29.

\bibitem{DSP19}
C. de Seguins Pazzis,
On the linear preservers of Schur matrix functionals,
{\it Linear Algebra Appl.}, {\bf 567}(2019) 63--117.

\bibitem{GM01}
B. Gelbord and R. Meshulam, Spaces of $p-$vectors of bounded rank,
{\it Israel J. of Math.}, {\bf 126}(2001) 129--139.	

\bibitem{Flanders62}
H. Flanders, On spaces of linear transformations with bounded rank,
{\it J. London Math. Soc.}, {\bf 37}(1962) 10--16.

\bibitem{GutSp} A. E. Guterman and I. A. Spiridonov, Permanent Polya problem for additive surjective maps,
{\it Linear Algebra Appl.}, {\bf 599}(2020) 140--155.

\bibitem{LQWWZ22}
Y. Li, Y. Qiao, A. Wigderson, Y. Wigderson and C. Zhang,
Connections between graphs and matrix spaces, arXiv:2206.04815.

\bibitem{Lovasz}
L. Lov\'{a}sz, {\it Combinatorial problems and exercises}, Corrected reprint of the 1993 second edition. AMS Chelsea Publishing, Providence, RI, 2007.

\bibitem{M85}
R. Meshulam, On the maximal rank in a subspace of matrices,
{\it Quarterly Journal of Mathematics Oxford}, {\bf 36}(1985) 225--229.

\bibitem{M17}
R. Meshulam,
Maximal rank in matrix spaces via graph matchings, {\it Linear Algebra and Appl.},
{\bf 529}(2017) 1--11.

\bibitem{West96}
D. B. West, Introduction to graph theory. Prentice Hall, Inc., Upper Saddle River, NJ, 1996.

\end{thebibliography}
\end{document}